\documentclass[reqno]{amsart}
\usepackage{amssymb}
\usepackage{hyperref}
\usepackage[all]{xy}

\newtheorem{thm}{Theorem}[section]
\newtheorem{prop}[thm]{Proposition}
\newtheorem{lem}[thm]{Lemma}
\newtheorem{cor}[thm]{Corollary}

\theoremstyle{definition}
\newtheorem{rem}[thm]{Remark}
\newtheorem{defn}[thm]{Definition}

\newtheorem{exas}[thm]{Examples}

\numberwithin{equation}{section}

\newfont{\cyrr}{wncyr10}

\def\Sig#1{\Sigma_{#1}}

\def\k{k}

\def\J{{\mathcal{J}}}
\def\I{{\mathcal{I}}}
\def\K{{\mathcal{K}}}
\def\cC{{\mathcal{C}}}
\def\O{{\mathcal{O}}}

\def\Z{\mathbb{Z}}
\def\Q{\mathbb{Q}}

\def\Gm{\mathbb{G}_m}
\def\Hom{\mathrm{Hom}}
\def\rk{\mathrm{rank}}

\def\Gal{\mathrm{Gal}}
\def\Res{\mathrm{Res}}

\def\GL{\mathrm{GL}}
\def\End{\mathrm{End}}
\def\Aut{\mathrm{Aut}}
\def\Spec{\mathrm{Spec}}

\def\coker{\mathrm{coker}}

\def\isom{\,\mbox{$\xrightarrow{\hskip6pt}\hskip-13pt\raisebox{3.5pt}{\footnotesize$\sim$}\hskip6pt$}}
\def\hookto{\hookrightarrow}
\def\map#1{\;\xrightarrow{#1}\;}
\def\dirsum#1{\underset{#1}{\textstyle\bigoplus}}

\def\onto{\twoheadrightarrow}
\def\s{\mathrm{s}}
\def\ks{\k^\s}
\def\bmu{{\boldsymbol{\mu}}}

\def\too{\longrightarrow}
\def\semidirect{\rtimes}

\def\mod#1{\text{\bf Mod}_\Z(#1/\k)}
\def\cag{\text{\bf CAG}}

\def\mod#1#2{\text{\bf FMod}_{#1}(#2)}
\def\modok{\mod{\O}{\k}}
\def\modlk{\mod{\O}{L/\k}}
\def\cag#1#2{\text{\bf CAG}_{#1}(#2)}
\def\cagok{\cag{\O}{\k}}
\def\Gm{\mathbb{G}_m}

\def\ItV{\I \otimes_\O V}
\def\JtV{\J \otimes_\O V}

\title{Twisting commutative algebraic groups}

\author{B.\ Mazur}
\address{Department of Mathematics,
Harvard University,
Cambridge, MA 02138 USA}
\email{mazur\char`\@math.harvard.edu}

\author{K.\ Rubin}
\address{Department of Mathematics,
University of California Irvine,
Irvine, CA 92697 USA}
\email{krubin\char`\@uci.edu}

\author{A.\ Silverberg}
\address{Department of Mathematics,
University of California Irvine,
Irvine, CA 92697 USA}
\email{asilverb\char`\@uci.edu}

\thanks{Mazur is supported by NSF grant DMS-0403374, Rubin by NSF grant DMS-0457481,
and Silverberg by NSA grant H98230-05-1-0044.}

\begin{document}

\begin{abstract}
If $V$ is a commutative algebraic group over a field $\k$, 
$\O$ is a commutative ring that acts on $V$, and $\I$ is a 
finitely generated free $\O$-module with a right action of 
the absolute Galois group of $\k$, then there is 
a commutative algebraic group $\ItV$ over $\k$, 
which is a twist of a power of $V$.  
These group varieties have applications to cryptography 
(in the cases of abelian varieties and algebraic tori over finite fields) 
and to the arithmetic of abelian varieties over number fields.
For purposes of such applications we devote this article to making 
explicit this tensor product construction and its basic properties.
\end{abstract}

\maketitle

\section*{Introduction}

In this paper we study twists of powers of commutative algebraic groups.  
We have been using and proving special cases of these results elsewhere, 
and believe that it would be useful to have a complete theory and 
complete proofs in the literature in one place.
Examples of applications of these twists 
that already appear in the literature 
include: to polarizations on abelian varieties \cite{howe}, 
to cryptography (\cite{freyfin,ssavs}
in the case of abelian varieties over finite fields and \cite{hwpaper} in
the case of algebraic tori over finite fields), 
to constructing abelian varieties over number fields with
Shafarevich-Tate groups of nonsquare order \cite{stein},
and to bounding below the Selmer rank of abelian varieties 
over dihedral extensions of number fields \cite{bigselmer}.

Suppose $V$ is a commutative algebraic group over a field $\k$, 
$\O$ is a commutative ring that acts on $V$, and $\I$ is a 
finitely generated free $\O$-module with a right action of 
the absolute Galois group $G_\k$ of $\k$.
We will define a commutative algebraic group over $\k$ that we will 
denote $\ItV$, which is a twist of a power of $V$.

The general theory underlying the construction of $\ItV$ is given 
in the standard sources discussing the homological algebra of 
tensor products of sheaves for the {\'e}tale topology  (see \cite{SGA4},
particularly Proposition 12.7 on p.~205).
In that language, $\ItV$ is a tensor 
product in the category of sheaves on the big \'etale site over $\Spec(\k)$.  
   This tensor product construction was introduced by Serre
   (\S2 of \cite{serreten}) in the case where $V$ is an elliptic curve
   with complex multiplication by $\O$ and $\I$ is a projective
   $\O$-module with trivial Galois action.  It is discussed in
   detail by Conrad (\S7 of \cite{conrad}) in the case where $V$ is a
   group scheme with $\O$-action and $\I$ is a projective
   $\O$-module with trivial Galois action.
Our main objective is to record, in some detail and in a usable way, 
the basic features regarding the operation of tensoring  an abelian 
group scheme over $\k$ endowed with a ring $\O$ of operators
(viewing the group scheme as a sheaf for the big {\'e}tale topology over $\Spec(\k$)) 
with a locally constant sheaf over $\Spec(k)$ of free $\O$-modules of finite rank, 
noting that the new sheaf given by this tensor product construction is 
again representable as a group scheme over $\k$ with an $\O$-action. 
To effectively make use of this construction in our applications, 
we found that we must pin things down very explicitly.  
For ease of reading we provide in this article a largely 
self-contained treatment.

In \S\ref{twisting} and \S\ref{twistprop}, 
following Milne \cite{milne} who dealt with the case of 
abelian varieties, we give a concrete definition 
of $\ItV$, and we prove some important properties.  
Some basic examples are given in Examples \ref{twistexs}.
Our construction is functorial in both $\I$ and $V$ (Theorem \ref{ithm}), 
and if $\I$ decomposes (up to finite index) as $\oplus_i \J_i$, 
then $\ItV$ is isogenous to $\oplus_i(\J_i \otimes V)$ 
(see Corollary \ref{decompcor}).  Theorem \ref{vvntl} describes 
the action of $G_\k$ on the torsion points of $\ItV$. 
We include a more general explicit construction of 
$\ItV$, without the assumption that $\I$ is a free $\O$-module, 
in an appendix.

If $L$ is a finite Galois extension of $\k$ and $G := \Gal(L/\k)$, 
then $\Z[G] \otimes_\Z V$ is the restriction of scalars $\Res^L_\k V$.
Theorem \ref{decprop} shows 
that $\Res^L_\k V$ is isogenous to $\oplus_\rho (\I_\rho \otimes_\Z V)$, 
where $\rho$ runs through the irreducible rational representations 
of $G$ and $\I_\rho$ is the intersection of $\Z[G]$ with 
the $\rho$-isotypic component of $\Q[G]$.  
In \S\ref{abelian} we restrict to the case where $L/\k$ is abelian, 
which is the case of interest in many of the applications.  
Similar results were obtained by Diem and Naumann \cite{diem} 
in the case of abelian varieties.  
In \S\ref{sdprodsect} we study cases where $\Gal(L/\k)$ is a
semi-direct product, which are needed for the applications in 
\cite{bigselmer}.
We study finite group actions on $\ItV$ in \S\ref{gpactsect};
these results have cryptographic significance in the case of
algebraic tori.

We thank Dick Gross for drawing our attention to 
Conrad's paper \cite{conrad}.
   
\subsection*{Notation}
Let $\Z^+$ denote the set of positive integers.  
If $\k$ is a field, $\ks$ will denote a separable closure of 
$\k$ and $G_\k := \Gal(\ks/\k)$.  
In this paper ``ring'' will always mean ring with identity, 
and ``commutative algebraic group'' will always
mean a commutative algebraic group variety (not
necessarily connected).

If $n \in \Z^+$, let $\bmu_n$ denote the group of 
$n$-th roots of unity in $\bar{\Q}$.
If $G$ is a finite group, then $\Z[G]$ will denote the group ring, 
except that $\Z[\bmu_n]$ denotes the ring of integers of the 
cyclotomic field $\Q(\bmu_n)$.

Suppose $\k$ is a field and $\O$ is a commutative ring.  
We consider two categories:

\begin{itemize}
\item
$\modok$ is the category whose objects are finitely 
generated free $\O$-modules with a continuous {\em right} action of $G_\k$, 
and whose morphisms are $G_\k$-equi\-vari\-ant $\O$-module homomorphisms 
(the modules are given the discrete topology, so a continuous $G_\k$-action 
is one that factors through a finite extension of $\k$).
\item
$\cagok$ is the category whose objects are 
commutative algebraic groups $V$ over $\k$ with an action of $\O$, i.e., 
a ring homomorphism $\O \to \End_\k(V)$, and whose morphisms are $\O$-equivariant 
homomorphisms defined over $\k$. 
\end{itemize}

If $\I,\J \in \modok$ we will view $\Hom_\O(\I,\J)$ as a left $G_k$-module, 
where for $f \in \Hom_\O(\I,\J)$, $\gamma \in G_k$, and $x \in \I$, we define 
$(f^\gamma)(x) = f(x\gamma)\gamma^{-1}$.  
View $\O \in \modok$ with trivial $G_k$-action.

\section{Twisting commutative algebraic groups}
\label{twisting}

Fix a field $\k$ and a commutative ring $\O$.  In this section we 
construct a functor $\modok \times \cagok \to \cagok$, 
which we will denote by $(\I,V) \mapsto \ItV$.  
This construction appears in \S2 of \cite{milne} when $V$ is an abelian variety.

\begin{defn}
\label{ItVdef}
Suppose $V \in \cagok$ and $\I \in \modok$.  Define the {\em $\I$-twist} 
$\ItV$ of $V$ as follows.
Let $d = \rk_\O(\I)$, and fix an $\O$-module isomorphism $j : \O^d \isom \I$. 
The homomorphism $\O \to \End_\k(V)$ induces  
$$
H^1(\k,\GL_d(\O)) \too H^1(\k,\Aut_\k(V^d)) \too H^1(\k,\Aut_{\ks}(V^d)),
$$
and we let $c_\I \in H^1(\k,\Aut_{\ks}(V^d))$ be the image of the 
cocycle $(\gamma \mapsto j^{-1} \circ j^\gamma)$ 
under this composition.  
Define $\ItV$ to be the twist of $V^d$ by the cocycle $c_\I$.  Namely, 
by Corollaire to Proposition 5 on p.~131 in \S{III-1.3} of \cite{serrecg} 
(see also \S3.1 of \cite{voskbk}), there is a pair 
$(\ItV,\phi)$ (unique up to isomorphism) where $\ItV \in \cagok$ and 
$\phi : V^d \isom \ItV$ is an isomorphism defined over $\ks$ such that for 
every $\gamma \in G_\k$, 
\begin{equation}
\label{phif}
c_\I(\gamma) = \phi^{-1} \circ \phi^\gamma.
\end{equation}
\end{defn}

\begin{rem}
\label{doL}
Suppose $L$ is a separable extension of $\k$ and $G_L$ acts trivially 
on $\I$.  Then $j^\gamma = j$ for all $\gamma \in G_L$, so $c_\I(\gamma) = 1$, 
so $\phi^\gamma = \phi$ by \eqref{phif}.  Thus 
the isomorphism $\phi : V^d \isom \ItV$ is defined over $L$. 
 
If we choose a different $\O$-module isomorphism $j' : \O^d \isom \I$ in 
Definition \ref{ItVdef}, then $j' = j \circ \alpha$ for some 
$\alpha \in \GL_d(\O)$.  The cocycles $\gamma \mapsto j^{-1} j^\gamma$ 
and $\gamma \mapsto (j')^{-1} (j')^\gamma = \alpha^{-1}j^{-1} j^\gamma \alpha^\gamma$ 
represent the same class 
in $H^1(\k,\GL_d(\O))$, so they give rise to the same class 
$c_\I \in H^1(\k,\Aut_{\ks}(V^d))$.  Thus $\ItV$ is independent of the choice of $j$.
\end{rem}

If $L/\k$ is a Galois extension,  
$V$ is a commutative algebraic group over $k$, $\I \in \modok$, 
and $A$ is a commutative $\k$-algebra, let $\gamma \in G_\k$ act on 
$A \otimes_\k L$ as $1 \otimes \gamma$ and on 
$\I \otimes_\O (V(A \otimes_\k L))$ as $\gamma^{-1} \otimes (1 \otimes \gamma)$.

\begin{lem}
\label{prerepble}
Suppose $\I \in \modok$, $V \in \cagok$, and 
$L$ is a Galois extension of $\k$ such that 
$G_L$ acts trivially on $\I$.  
Fix an $\O$-module isomorphism $j : \O^d \isom \I$, and let 
$\phi : V^d \isom \ItV$ be as in Definition \ref{ItVdef}. 
Then for every commutative $\k$-algebra $A$,
the composition 
$$
(\ItV)(A \otimes_\k L) \isom \I \otimes_\O (V(A \otimes_\k L))
$$
of the sequence of $\O$-module isomorphisms
$$
(\ItV)(A \otimes_\k L)
    \map{\phi^{-1}} V^d(A \otimes_\k L) 
    \isom \O^d \otimes_\O (V(A \otimes_\k L)) 
    \map{j \otimes 1} \I \otimes_\O (V(A \otimes_\k L))
$$
is a $G_\k$-equivariant $\O$-module isomorphism 
that is independent of $j$ and is functorial in $A$, $V$, $\I$, and $L$.
\end{lem}

\begin{proof}
Remark \ref{doL} shows that $\phi$ is defined over $L$, and 
therefore $\phi(V^d(A \otimes_\k L))=(\ItV)(A \otimes_\k L)$.
The $G_\k$-equivariance of the composition and
the independence of $j$ follow from \eqref{phif} and the 
definition of $c_\I$.  
The functoriality is clear.
\end{proof}

\begin{thm}
\label{repble}
Suppose $\I \in \modok$ and $V \in \cagok$.  
Let $L$ be a Galois extension of $\k$ such that 
$G_L$ acts trivially on $\I$.   
Then $\ItV$ represents the functor on commutative $\k$-algebras 
$A \mapsto (\I \otimes_\O (V(A \otimes_\k L)))^{\Gal(L/\k)}$.
More precisely, for every commutative $\k$-algebra $A$,
the isomorphism of Lemma \ref{prerepble} restricts to a
functorial group isomorphism
$$
(\ItV)(A)  \cong  (\I \otimes_\O (V(A \otimes_\k L)))^{\Gal(L/\k)}.
$$
\end{thm}

\begin{proof}
This follows directly from Lemma \ref{prerepble}, since 
$(A \otimes_\k L)^{G_\k} = A$ and $G_L$ acts trivially on $\I$ and $L$.
\end{proof}

\begin{exas}
\label{twistexs}
\begin{enumerate}
\item
Suppose $0 \le d \in \Z$, and $\I = \O^d$ with trivial Galois action.  Then 
$\ItV = V^d$.
\item
Suppose $\chi$ is a quadratic character of $G_\k$, and 
$\I$ is a free, rank-one $\Z$-module with $G_\k$ acting via $\chi$.  
Then $\I \otimes_\Z V$ is the quadratic twist of $V$ by $\chi$.
More generally, if $\O = \Z[\bmu_n]$, $\chi : G_\k \to \bmu_n$ 
is a homomorphism, and $\I$ 
is a free, rank-one $\O$-module with $G_\k$ acting via $\chi$, then 
$\ItV$ is the twist of $V$ by $\chi^{-1}$ 
(in this case the cocycle $c_\I$ is $\chi^{-1}$).
\item
Suppose $V = \Gm$, the multiplicative group, and $\I$ is a free $\Z$-module.  
Then $\I \otimes_\Z V$ is the algebraic torus whose 
character module $\Hom(\I \otimes_\Z V, \Gm)$ is $\Hom(\I,\Z)$.  
See Corollary \ref{cmcor} below, and Example 6 in \S3.4 of \cite{voskbk}.
\item
If $L/k$ is a finite Galois extension then 
$\O[\Gal(L/\k)] \otimes_\O V = \Res^L_\k V$ 
(see Proposition \ref{rescpr1} below).  
\end{enumerate}
\end{exas}

\begin{prop}
\label{newjdiag}
Suppose $\I,\J \in \modok$ and $V,W \in \cagok$.
\begin{enumerate}
\item
There is a functorial $G_\k$-equivariant $\O$-module isomorphism
$$
\Hom_\O(\I,\J) \otimes_\O \Hom_{\ks}(V,W) 
    \isom \Hom_{\ks}(\ItV, \J \otimes_\O W).
$$
\item
The isomorphism of (i) restricts to an injective homomorphism 
$$
\Hom_{\O[G_k]}(\I,\J) \otimes_\O \Hom_{\k}(V,W) 
    \hookto \Hom_{\k}(\ItV, \J \otimes_\O W).
$$
\end{enumerate}
\end{prop}

\begin{proof}
Fix $\O$-module isomorphisms $\O^n \cong \I$ and $\O^m \cong \J$.  
These isomorphisms induce (see Definition \ref{ItVdef}) 
isomorphisms $V^n \isom \ItV$ and $W^m \isom \J \otimes_\O W$ 
defined over $\ks$, which induce isomorphisms 
\begin{multline*}
\Hom_\O(\I,\J) \otimes_\O \Hom_{\ks}(V,W) 
    \isom M_{m \times n}(\O) \otimes_\O \Hom_{\ks}(V,W) \\
    \isom M_{m \times n}(\Hom_{\ks}(V,W)) 
    \isom \Hom_{\ks}(V^n,W^m) 
    \isom \Hom_{\ks}(\ItV, \J \otimes_\O W).
\end{multline*} 
The proof of $G_\k$-equivariance is similar to the proof of 
$G_\k$-equivariance in Lemma \ref{prerepble}.
This proves (i), and (ii) follows since $G_\k$ acts trivially on
$
\Hom_{\O[G_k]}(\I,\J) \otimes_\O \Hom_{\k}(V,W).
$
\end{proof}

\begin{cor}
\label{jdrem}
Suppose $\I,\J \in \modok$ and $V \in \cagok$.
\begin{enumerate}
\item
The isomorphism of Proposition \ref{newjdiag}(i) with 
$W = V$ and the identity map in $\Hom_\k(V,W)$ gives a 
functorial $G_\k$-equivariant $\O$-module homomorphism
$$
\Hom_\O(\I,\J) \too \Hom_{\ks}(\ItV, \J \otimes_\O V).
$$
\item
The map of (i) restricts to a homomorphism
$$
\Hom_{\O[G_\k]}(\I,\J) \too \Hom_{\k}(\ItV, \J \otimes_\O V).
$$
\item
If the map $\O \to \End_\k(V)$ is injective, then the maps 
in (i) and (ii) are injective.
\end{enumerate}
\end{cor}

\begin{proof}
Assertions (i) and (ii) follow directly from Proposition \ref{newjdiag}.  
For (iii), tensoring the injection $\O \hookto \End_\k(V)$ with the free 
$\O$-module $\Hom_\O(\I,\J)$ shows that 
$\Hom_\O(\I,\J) \hookto \Hom_\O(\I,\J) \otimes_\O \End_\k(V)$ is injective.  
Now (iii) follows from the injectivity in Proposition \ref{newjdiag}.
\end{proof}

If $f \in \Hom_\O(\I,\J)$ we will often write $f_V$ for the image of $f$ 
under the map of Corollary \ref{jdrem}(i).

\begin{thm}
\label{ithm}
The map $(\I,V) \mapsto \ItV$ is a functor from 
$\modok \times \cagok$ to $\cagok$.
\end{thm}

\begin{proof}
This follows directly from Proposition \ref{newjdiag}(ii).
\end{proof}

\begin{cor}
\label{obvcor}
Suppose $\I, \J \in \modok$, $V \in \cagok$, and 
$\k \subseteq F \subseteq \ks$.  If $\I$ and $\J$ 
are isomorphic as $\O[G_F]$-modules, then the group varieties $\ItV$ 
and $\J \otimes_\O V$ are isomorphic over $F$.
\end{cor}

\begin{proof}
If $f : \I \to \J$ is a $G_F$-equivariant isomorphism, then the 
image of $f$ under the functorial map 
of Corollary \ref{jdrem}(i) 
is an isomorphism over $F$ from $\ItV$ to $\JtV$.
\end{proof}

\begin{cor}
\label{cmcor}
If $\I \in \mod{\Z}{\k}$, then 
$\Hom_{\ks}(\I \otimes_\Z \Gm, \Gm) \cong \Hom_\Z(\I,\Z)$.
\end{cor}

\begin{proof}
Apply Proposition \ref{newjdiag}(i) with $\J = \O = \Z$ and $V = W = \Gm$.
\end{proof}

\section{Properties of the twists $\ItV$}
\label{twistprop}

For this section, fix a field $\k$, a commutative ring $\O$, 
and a commutative algebraic group $V \in \cagok$.

\begin{thm}
\label{easycor}
Suppose $\I \in \modok$.  Then:
\begin{enumerate}
\item
$\ItV$ is a commutative algebraic group of dimension $\rk_\O(\I)\dim(V)$, 
\item 
$\ItV$ is connected if and only if $V$ is connected., 
\item
if $L$ is a separable extension 
of $\k$ and $G_L$ acts trivially on $\I$, then $\ItV$ is 
isomorphic over $L$ to $V^{\rk_\O(\I)}$.
\end{enumerate}
\end{thm}

\begin{proof}
Fix a separable extension $L/\k$ such that $G_L$ acts trivially on $\I$.  
Since $\I$ is isomorphic as a $G_L$-module to $\O^{\rk_\O(\I)}$, 
$\ItV$ is isomorphic over $L$ to 
$\O^{\rk_\O(\I)} \otimes_\O V = V^{\rk_\O(\I)}$ by Corollary \ref{obvcor}, 
giving (iii).  The remaining assertions follow easily.
\end{proof}

Suppose $n \in \Z^+$. If $B$ is an abelian group, let $B[n]$ denote the 
subgroup of elements of order dividing $n$ in $B$.
If $W$ is a commutative algebraic group over $\k$, 
let $W[n]$ denote the $G_\k$-module $W(\ks)[n]$, and if $\ell$ is a prime let 
$$T_\ell(W) := \varprojlim_m W[\ell^m],$$
the $\ell$-adic Tate module of $W$.

\begin{thm}
\label{vvntl}
Suppose $\I \in \modok$, $n \in \Z^+$, and $\ell$ is prime.  
Then there are $G_\k$-equivariant isomorphisms (with $\gamma \in G_\k$ 
acting on the right-hand sides as $\gamma^{-1} \otimes \gamma$), 
functorial in $\I$ and $V$,
\begin{enumerate}
\item
$(\ItV)(\ks) \cong \I \otimes_\O (V(\ks))$,
\item
$(\ItV)[n] \cong \I \otimes_\O (V[n])$,
\item
$T_\ell(\ItV) \cong \I \otimes_\O (T_\ell(V))$.
\end{enumerate}
\end{thm}

\begin{proof}
(See Proposition 6(b) of \cite{milne}.)
The first assertion follows from Lemma \ref{prerepble} with $A = \k$ 
and $L = \ks$.
Since $\I$ is a free $\O$-module, 
$$
(\I \otimes_\O (V(\ks)))[n] \cong \I \otimes_\O (V[n]),
$$
so (ii) follows from (i), and (iii) follows by taking the inverse limit
of (ii) with $n=\ell^m$.
\end{proof}

\begin{lem}
\label{twex}
Suppose $\I, \J \in \modok$ with $\I \subseteq \J$ and $\J/\I$ is 
free (as an $\O$-module).  
Then the induced sequence of commutative 
algebraic groups over $\k$
$$
0 \too \ItV \too \J \otimes_\O V \too (\J/\I) \otimes_\O V \too 0
$$
is exact.
\end{lem}

\begin{proof}
Since $\J/\I$ is free, there is an $\O$-module 
isomorphism $\J \cong \I \oplus (\J/\I)$.  
It follows by Corollary \ref{obvcor} 
that $\JtV \cong (\ItV) \oplus ((\J/\I) \otimes_\O V)$ over $\ks$, so the 
sequence of the lemma is a (split) exact sequence over $\ks$.  
But then the sequence is exact over $\k$.
\end{proof}

Define a {\em $\k$-isogeny} in $\cagok$ or in $\modok$ to be a
$\k$-morphism whose kernel and cokernel are annihilated by some 
positive integer.

\begin{lem}
\label{isoglem}
If $\I, \J \in \modok$ and $s : \I \to \J$ is a $\k$-isogeny,
then the induced map $s_V : \ItV \to \JtV$ is a $\k$-isogeny.
\end{lem}

\begin{proof}
Suppose $n \in \Z^+$ 
is such that $n \cdot \ker(s) = 0$ and $n \cdot \coker(s) = 0$.  
Then there is a $\k$-isogeny 
$t : \J \to \I$ such that $t \circ s$ and $s \circ t$ are both 
multiplication by $n^2$, so 
$s_V \circ t_V \in \End_\k(\J \otimes_\O V)$ and 
$t_V \circ s_V \in \End_\k(\ItV)$ 
are both multiplication by $n^2$.  Therefore $s_V$ is a $\k$-isogeny.
\end{proof}

\begin{cor}
\label{decompcor}
Suppose $\I, \J_1, \ldots, \J_t \in \modok$, and 
$\I \otimes_\Z \Q \cong \oplus_{i=1}^t (\J_i \otimes_\Z \Q)$ 
as $\O[G_\k]$-modules.  Then 
$\ItV$ is $\k$-isogenous to $\oplus_{i=1}^t (\J_i \otimes_\O V)$.
\end{cor}

\begin{proof}
In this case $\I$ is $\k$-isogenous to $\oplus_i \J_i$, so by 
Lemma \ref{isoglem}, $\ItV$ is $\k$-isogenous to 
$(\oplus_i \J_i) \otimes_\O V \cong \oplus_i (\J_i \otimes_\O V)$.
\end{proof}

\begin{prop}
\label{repeat}
Suppose $\I, \J \in \modok$.  
Then there is a natural isomorphism 
$(\I \otimes_\O \J) \otimes_\O V \cong \I \otimes_\O (\J \otimes_\O V)$
over $\k$.
\end{prop}

\begin{proof}
Suppose $A$ is a commutative $\k$-algebra.  Then applying 
Theorem \ref{repble} and Lemma \ref{prerepble} with $L = \ks$ 
(suppressing the subscripts $\O$ and $\k$ from the tensor products)
\begin{align*}
(\I \otimes (\J \otimes V))(A) 
    & \cong (\I \otimes ((\J \otimes V)(A \otimes \ks)))^{G_\k} \\
    & \cong (\I \otimes (\J \otimes (V(A \otimes \ks)))^{G_\k} \\
    & = ((\I \otimes \J) \otimes (V(A \otimes \ks)))^{G_\k} \\
    & \cong ((\I \otimes \J) \otimes V)(A).
\end{align*}
These isomorphisms are functorial in $A$, so the proposition follows 
from a variant of the Yoneda Lemma (see for example Proposition VI-2
of \cite{eisenbud}).
\end{proof}

\section{Annihilator modules}
\label{lemssect}

The results of this section will be used in 
\S\ref{resn} and \S\ref{gpactsect}.

Fix a finite Galois extension $L/\k$, a commutative ring $\O$, and 
a commutative algebraic group $V \in \cagok$, and let $G := \Gal(L/\k)$.  
Let $\modlk$ denote the full subcategory of $\modok$ whose objects 
are the $\O[G_\k]$-modules in $\modok$ on which $G_L$ acts trivially.

\begin{defn}
\label{hatdef}
For $\I \in \modlk$, define a left $\O[G]$-module
$$
\hat\I := \Hom_{\O[G]}(\I,\O[G])
$$  
with $\O[G]$ acting by 
$(\alpha \cdot f)(x) = \alpha \cdot f(x)$ 
for every $\alpha \in \O[G]$, $f \in \hat{\I}$, and $x \in \I$.
Also define an $\O$-module homomorphism 
$\pi : \O[G] \to \O$ by $\pi(\sum_{g\in G} a_g g) = a_1$.
\end{defn}

Part (i) of the following lemma shows that $\hat\I$ 
is independent of the choice of $L$.

\begin{lem}
\label{hatexact}
\begin{enumerate}
\item
For each $\K \in \modlk$, the map  
$f \mapsto \pi \circ f$ defines an isomorphism of left $G_k$-modules
$$
\hat\K \isom \Hom_\O(\K,\O).
$$
\item
If $\I, \J \in \modlk$, $\I \subseteq \J$, and $\J/\I$ is a projective   
$\O$-module, then the canonical sequence
$
0 \to \widehat{\J/\I} \to \hat\J \to \hat\I \to 0
$
is exact.
\end{enumerate}
\end{lem}

\begin{proof}
For (i), see for example Proposition VI.3.4 of \cite{brown}.
Assertion (ii) follows from (i), 
since the exact sequence  of $\O$-modules $0 \to \I \to \J \to \J/\I \to 0$ 
splits if $\J/\I$ is projective.
\end{proof}

\begin{lem}
\label{iperp}
Suppose $\I, \J \in \modlk$, $\I \subseteq \J$, and $\J/\I$ is a free   
$\O$-module.  
Then 
$$
\ItV = \bigcap_{f \in \widehat{\J/\I}} \ker(f_V : \JtV \to \O[G] \otimes_\O V),
$$
where $f_V$ is the image of $f$ under the map 
$$
\Hom_{\O[G]}(\J/\I,\O[G]) \hookto \Hom_{\O[G]}(\J,\O[G]) 
    \too \Hom_\k(\JtV, \O[G] \otimes_\O V)
$$
coming from Corollary \ref{jdrem}(ii).
\end{lem}

\begin{proof}
Choose an $\O$-basis $\{f_1,\ldots,f_d\}$ of $\Hom_\O(\J/\I,\O)$.  For 
every $i$ let $\phi_i \in \widehat{\J/\I}$ be the 
inverse image of $f_i$ under the isomorphism of 
Lemma \ref{hatexact}(i) (with $\K = \J/\I$), so $\pi \circ \phi_i = f_i$, 
with $\pi$ defined in Definition \ref{hatdef}.  
Then there is a commutative diagram of $\O$-modules, with the top line 
an exact sequence of $\O[G]$-modules
\begin{equation}
\label{sek}
\parbox{4in}{\xymatrix{
0 \ar[r] & \J/\I \ar^-{\oplus \phi_i}[r] \ar^(.575){\hskip -5pt \oplus{f_i}}_{\cong}[dr] 
    & \O[G]^d \ar[r] \ar^-(.4){\pi^d}[d] & \cC \ar[r] & 0. \\
&& \O^d
}}
\end{equation}
Then $(\oplus f_i)^{-1} \circ \pi^d$ gives a splitting of the top exact 
sequence, so $\cC$ is isomorphic as an $\O$-module to the kernel of $\pi^d$, 
which is free.  
Therefore we can apply Lemma \ref{twex} both to the top line of \eqref{sek} 
and to the exact sequence $0 \to \I \to \J \to \J/\I \to 0$ 
to obtain an exact sequence
$$
0 \to \ItV \to \JtV \map{\oplus (\phi_i)_V} 
    (\O[G] \otimes_\O V)^d \to \cC \otimes_\O V \to 0.
$$
By Lemma \ref{hatexact}(i), $\phi_1,\ldots,\phi_d$ 
generate $\widehat{\J/\I}$, so 
$$
\ItV = \ker(\oplus(\phi_i)_V) = \cap_{f \in \widehat{\J/\I}} \ker(f_V).
$$ 
\end{proof}
\begin{defn}
If $\I$ is a right ideal of $\O[G]$, let $\I^\perp$ denote the 
left annihilator of $\I$, i.e., $\I^\perp$ is the left ideal of $\O[G]$
defined by
$$
\I^\perp := \{\alpha \in \O[G] : \alpha \I = 0\}.
$$
A (right or left) ideal $\I$ of $\O[G]$ is 
{\em saturated} if $\O[G]/\I$ is a projective $\O$-module.
\end{defn}

A finitely generated 
$\Z$-module is projective (or equivalently, free) 
if and only if it is torsion-free.  Thus when $\O = \Z$, 
intersecting with $\Z[G]$ (inversely, tensoring with $\Q$) gives a one-to-one 
correspondence between the ideals of $\Q[G]$ and  the saturated ideals of $\Z[G]$.

\begin{lem}
\label{endl}
Let $\lambda : \O[G] \to \widehat{\O[G]}$ be the ring isomorphism 
that sends $\alpha \in \O[G]$ to left multiplication by $\alpha$.  Then:
\begin{enumerate}
\item
If $\I$ is a right ideal of $\O[G]$ then the restriction of
$\lambda$ induces an isomorphism 
$$
\I^\perp \isom \widehat{\O[G]/\I}.
$$
\item
If $\I$ is a saturated right ideal of $\O[G]$, 
then $\I = \{\alpha \in \O[G] : \I^\perp \cdot \alpha = 0\}$.
\item
If $\I$ is a saturated two-sided ideal of $\O[G]$, then $\lambda$ induces 
an isomorphism $\O[G]/\I^\perp \isom \End_{\O[G]}(\I)$ 
(and $\End_{\O[G]}(\I) = \hat\I$).
\end{enumerate}
\end{lem}

\begin{proof}
Suppose $\I$ is a right ideal of $\O[G]$.
The map $\widehat{\O[G]} \to \O[G]$ defined by $f \mapsto f(1)$
is a right and left inverse of $\lambda$.  Thus $\lambda$
is an isomorphism and there is a commutative diagram with exact rows
\begin{equation}
\label{ex?}
\parbox{5in}{\xymatrix@R=20pt{
0 \ar[r] & \I^\perp \ar[r]\ar_{\lambda}[d]
    & \O[G] \ar[r]\ar_{\lambda}^{\cong}[d]
    & \O[G]/\I^\perp_{\phantom{_X}} \ar[r]\ar@{^(->}_{\lambda}[d] & 0 \\
0 \ar[r] & \widehat{\O[G]/\I} \ar[r] & \widehat{\O[G]} \ar[r] & \hat\I
}}
\end{equation}
where the right-hand vertical map is injective by definition of
$\I^\perp$.  The snake lemma 
shows that the left-hand vertical map is an isomorphism, which proves (i).

Now suppose $\I$ is saturated.
If $\beta \in \O[G] - \I$, then there is an $\O$-module 
homomorphism from $\O[G]/\I$ to $\O$ that is nonzero on $\beta$.
Now (ii) follows from (i), along with Lemma \ref{hatexact}(i) with
$\K = \O[G]/\I$.

Since $\I$ is saturated, by Lemma \ref{hatexact}(ii) the bottom
right-hand map of \eqref{ex?} is surjective, and hence the right-hand
vertical map is an isomorphism.
If $\I$ is a two-sided ideal, then
$\lambda(\O[G]/\I^\perp) \subseteq \End_{\O[G]}(\I) \subseteq \hat\I$,
so equality must hold and the proof of (iii) is complete.
\end{proof}

\section{Decomposing the restriction of scalars}
\label{resn}

In this section we decompose the restriction of scalars of a commutative 
algebraic group.  
Theorems \ref{decprop}, \ref{decomp}, and \ref{abprim0} 
were proved by Diem and Naumann 
in \S3.4 and \S3.5 of \cite{diem} in the case of abelian varieties.

Fix a finite Galois extension $L/\k$, a commutative ring $\O$, and 
a commutative algebraic group $V \in \cagok$, and let $G := \Gal(L/\k)$.  
Let $\Res^L_\k V$ denote the Weil restriction of scalars of $V$ from $L$ 
to $\k$ (see for example \S1.3 of \cite{weil} or \S 3.12 of \cite{voskbk}). 
Then for every commutative $\k$-algebra $A$ there is an isomorphism, 
functorial in $A$,
\begin{equation}
\label{resscal}
(\Res^L_\k V)(A) \cong V(A \otimes_\k L).
\end{equation} 

\begin{prop}
\label{rescpr1}
$\O[G] \otimes_\O V \cong \Res^L_\k V$ over $\k$.
\end{prop}

\begin{proof}
Let $\O^G := \oplus_{g \in G} \O$, $V^G := \oplus_{g \in G} V$, and 
for $g \in G$ let $p_g : V^G \to V$ be the projection onto the $g$ component.  
Using the $\O$-module isomorphism $j : \O^G \to \O[G]$ defined by 
$j((x_g)) = \sum_g x_g g^{-1}$, 
Definition \ref{ItVdef} gives a pair $(\O[G] \otimes_\O V, \phi)$ 
where 
$\phi : V^G \isom \O[G] \otimes_\O V$ is an isomorphism over $L$.  
Let $\eta := p_1 \circ \phi^{-1} : \O[G] \otimes_\O V \to V$.  
The cocycle $c_{\O[G]} \in H^1(\k,\Aut_{\ks}(V^G))$ of Definition \ref{ItVdef} 
satisfies $p_h \circ c_{\O[G]}(g) = p_{g^{-1}h}$ for every $g, h \in G$, 
so (using \eqref{phif}), 
$$
\eta^g = p_1 \circ (\phi^{-1})^g = p_1 \circ c_{\O[G]}(g)^{-1} \circ \phi^{-1} 
    = p_{g} \circ \phi^{-1}.
$$
Therefore $\oplus \eta^g :  \O[G] \otimes_\O V \to V^G$ is an isomorphism 
(it's equal to $\phi^{-1}$), so by the definition of $\Res^L_\k V$ 
in \S1.3 of \cite{weil}, 
$\O[G] \otimes_\O V \cong \Res^L_\k V$ over $\k$.
\end{proof}

For the rest of this section we will take $\O = \Z$ and write 
simply ``$\otimes$'' in place of ``$\otimes_\Z$''.  
The functorial map of Corollary \ref{jdrem}(ii) (with $\I=\J=\Z[G]$) 
and Proposition \ref{rescpr1} give natural ring homomorphisms 
\begin{equation}
\label{rescpr2}
\Z[G]  \cong  \End_{\Z[G]}(\Z[G]) \to \End_\k(\Z[G] \otimes V) 
    \isom \End_\k(\Res^L_\k V).
\end{equation}
If $\alpha \in \Z[G]$, then we denote its image under \eqref{rescpr2} 
by $\alpha_V \in \End_\k(\Res^L_\k V)$.

\begin{prop}
\label{altthm}
If $\I$ is a saturated right ideal of $\Z[G]$, then:
\begin{enumerate}
\item
$\I \otimes V = \bigcap_{\alpha \in \I^\perp} \ker(\alpha_V : \Res^L_\k V \to \Res^L_\k V)$.
\item
For every commutative $\k$-algebra $A$ there is a functorial isomorphism
$$
(\I \otimes V)(A) \cong \{v \in V(A \otimes_\k L) : \I^\perp \cdot v = 0\}.
$$
\item
If further $\I$ is a two-sided ideal, then there is a natural injective 
ring homomorphism 
$$
(\Z[G]/\I^\perp) \otimes \End_\k(V) \hookto \End_\k(\I \otimes V).
$$
\end{enumerate}
\end{prop}

\begin{proof}
By Lemma \ref{endl}(i), $\widehat{\Z[G]/\I} \cong \I^\perp$.
By Lemma \ref{iperp} (with $\J = \Z[G]$) and Proposition \ref{rescpr1}, we have (i).  Assertion (ii) 
follows from (i) and \eqref{resscal}. 

If $\I$ is a two-sided ideal, then
$\Z[G]/\I^\perp \isom \End_{\Z[G]}(\I)$ by Lemma \ref{endl}(iii).  
Now (iii) follows from Proposition \ref{newjdiag}(ii) 
(with $\J = \I$ and $W = V$).
\end{proof}

The group ring $\Q[G]$ is semisimple, and decomposes into a direct sum 
of minimal two-sided ideals 
\begin{equation}
\label{ssd}
\Q[G] = \dirsum{\rho} \Q[G]_\rho
\end{equation}
indexed by the irreducible rational representations $\rho$ of $G$.  
Here $\Q[G]_\rho$ is the $\rho$-isotypic component of $\Q[G]$, i.e., the 
sum of all left ideals of $\Q[G]$ isomorphic to $\rho$.

\begin{defn}
\label{twdef}
If $\rho$ is an irreducible finite-dimensional 
rational representation of $G_\k$,  
choose a finite Galois extension $L/\k$ such that 
$\rho$ factors through $G := \Gal(L/\k)$,  
define 
$$
\I_\rho := \Q[G]_\rho \cap \Z[G] \in \modok,
$$ 
and define the {\em $\rho$-twist of $V$} by
$$
V_\rho := \I_\rho \otimes V.
$$
\end{defn}

\begin{rem}
\label{rhoind}
Note that $\I_\rho$ is well-defined up to $\Z[G_\k]$-isomorphism, 
independent of the choice of $L$, and therefore $V_\rho$ is well-defined 
up to isomorphism over $\k$.  
Since $\Q[G]_\rho$ is a $\Q$-vector space, $\I_\rho$ is a saturated 
ideal of $\Z[G]$.
\end{rem}

\begin{thm}
\label{decprop}
Suppose $L/\k$ is a finite Galois extension, 
$V$ is a commutative algebraic group over $\k$, and
$G := \Gal(L/\k)$.
Then $\Res^L_\k V$ is isogenous over $\k$ to $\oplus_\rho V_{\rho}$, 
direct sum over all irreducible rational representations of $G$.
\end{thm}

\begin{proof}
This follows from Proposition \ref{rescpr1}, \eqref{ssd}, and Corollary \ref{decompcor} with 
$\I = \Z[G]$ and with
$
\{\J_1,\ldots,\J_t\} = \{\I_\rho : \text{$\rho$ an irreducible 
rational representation of $G$}\}.
$
\end{proof}

\section{Abelian twists}
\label{abelian}

Fix a finite abelian extension $L/\k$ and a commutative 
algebraic group $V$ over $\k$, and let $G = \Gal(L/\k)$ and $\O = \Z$.

The irreducible rational representations of $G$ are in one-to-one 
correspondence with the cyclic extensions of $\k$ in $L$.  
(See for example exercise 13.1 of \cite{serre}.)
Namely, if $\rho$ is an irreducible rational representation let $F_\rho$ be 
the fixed field of the kernel of $\rho$, and if $F$ is a cyclic extension 
of $\k$ in $L$ let $\rho_F$ (or $\rho_{F/\k}$, if we 
need to specify the field $\k$) denote the unique irreducible rational 
representation of $G$ with kernel $\Gal(L/F)$.  If $[F:\k] = d$ then 
$\dim{\rho_F} = \varphi(d)$, where $\varphi$ is the Euler $\varphi$-function.  

\begin{defn}
\label{vldef}
Suppose $F$ is a cyclic extension of $\k$ in $L$, and $\rho_F$ 
is the corresponding irreducible rational representation of $G$.  
Let $\Q[G]_F$ denote the $\rho_F$-isotypic component of $\Q[G]$, 
and let 
$$
\I_F := \Q[G]_F \cap \Z[G], \qquad V_F := \I_F \otimes V
$$
(these were denoted $\Q[G]_{\rho_F}$, $\I_{\rho_F}$, and 
$V_{\rho_F}$ in \eqref{ssd} and Definition \ref{twdef}).  
When necessary to specify the ground field $\k$, we will write 
$\I_{F/\k}$ and $V_{F/\k}$.  Let $R_F$ denote the 
maximal order of the field $\Q[G]_F$.
\end{defn}

By Remark \ref{rhoind}, $\I_F$ and 
$V_F$ are well-defined up to isomorphism, independent of the choice 
of field $L$ containing $F$, 
and $\I_F$ is saturated in $\Z[G]$.

The following result is a special case of Theorem \ref{decprop}.

\begin{thm}
\label{decomp}
If $L/\k$ is a finite abelian extension and $V$ is a commutative 
algebraic group over $\k$, then $\Res^L_\k V$ is isogenous over 
$k$ to $\oplus_F V_F$, direct sum over all 
cyclic extensions $F$ of $\k$ in $L$.
\end{thm}

If $\k \subseteq F \subseteq L$, let 
$$
N_{L/F} := \sum_{g \in \Gal(L/F)} g \in \Z[G].
$$
Define  
$$
\Omega_{L} := \{\text{fields $F$} : \k \subseteq F \subsetneq L\} \quad\supseteq\quad
\Omega'_{L} := \{F : \text{$\k \subseteq F \subsetneq L$, $[L:F]$ prime}\}.
$$
Then every element of $\Omega_L$ is a subfield of some element of $\Omega'_L$, 
and we define 
$$
W_L := \bigcap_{F \in \Omega_L} \ker (N_{L/F,V}) 
= \bigcap_{F \in \Omega'_L}\ker(N_{L/F,V})\subseteq \Res^L_\k V,
$$
where $N_{L/F,V} \in \End_k(\Res^L_\k V)$ is the image of $N_{L/F}$ 
under \eqref{rescpr2}.
We will see in Theorem \ref{abprim}(i) below that if $L/\k$ is cyclic, 
then $W_L = V_L$.  In the non-cyclic case we have the following.

\begin{prop}
If $L/\k$ is abelian but not cyclic, then $\dim(W_L) = 0$.
\end{prop}

\begin{proof}
Since $L/\k$ is not cyclic, there are a prime $p$ and a field 
$M$ such that $\k \subseteq M \subset L$ and $\Gal(L/M) \cong (\Z/p\Z)^2$.  
Since there are exactly $p+1$ degree $p$ extensions of $M$ in $L$, 
in $\Z[G]$ we have the identity
$$
\sum_{M \subsetneq F \subsetneq L} N_{L/F} = p + N_{L/M}.
$$
Since $W_L$ is in the kernel of all the norm maps in this identity, 
it follows that $W_L$ is contained 
in the kernel of multiplication by $p$, so $\dim(W_L) = 0$. 
\end{proof}

Suppose for the rest of this section that $L/\k$ is cyclic.  
Theorems \ref{abprim0},  \ref{abprim}, and \ref{abprim3} below 
are our main results about $V_L$ in the cyclic case.   
Let $r := |G| = [L:\k]$, and fix a generator $\tau$ of $G$.  
For $d \in \Z^+$ let $\Phi_d \in \Z[x]$ denote the $d$-th 
cyclotomic polynomial, and let $\Psi_d(x) := (x^d-1)/\Phi_d(x) \in \Z[x]$.  

\begin{lem}
\label{newlem}
\begin{enumerate}
\item
$\I_L= \Psi_r(\tau)\Z[G]$ and $\I_L^\perp=\Phi_r(\tau)\Z[G]$.
\item
Every isomorphism $\chi : G \isom \bmu_r$ induces a ring isomorphism $R_L \isom \Z[\bmu_r]$.  
This ring isomorphism is 
$G$-equivariant, with $g \in G$ acting on $\Z[\bmu_r]$ as multiplication 
by $\chi(g)$.
\item
The projection $\Q[G] \onto \Q[G]_L$ given by \eqref{ssd} 
induces a $G$-module isomorphism
$\Z[G]/\I_L^\perp \isom R_L$.
\item
$\I_L = \textstyle\prod_{\text{primes $\ell\mid r$}}(\zeta_\ell - 1) R_L$,
where for each prime $\ell$ dividing $r$, $\zeta_\ell$ is a
primitive $\ell$-th root of unity in $R_L$.
\end{enumerate}
\end{lem}

\begin{proof}
Let $S=\Q[x]/(x^r-1)\Q[x]$.
Since $G$ is cyclic of order $r$, the homomorphism $\eta : \Q[G] \to S$ 
that takes $\tau$ to $x$
is a $\Q[G]$-module isomorphism, where $\tau$ acts on 
$S$ as multiplication by $x$.  Since
$\Q[G]_L \cong \Q[x]/\Phi_r(x)\Q[x] \cong \Psi_r(x)S
    \subseteq S$
as $\Q[G]$-modules,
and $\Q[G]$ (and hence $S$) has a unique $\Q[G]$-submodule isomorphic to $\Q[G]_L$, 
we have $\eta(\Q[G]_L) = \Psi_r(x)S$.  
It follows that the isomorphism  $\eta : \Z[G] \cong \Z[x]/(x^r-1)\Z[x]$ 
maps $\I_L$ (resp., $\I_L^\perp$) 
isomorphically onto the ideal generated by $\Psi_r(x)$ (resp., by $\Phi_r(x)$).
Both assertions of (i) now follow.

If $\chi : G \isom \bmu_r$ is an isomorphism, then 
$\tau \mapsto x \mapsto \chi(\tau)$ 
induces isomorphisms $\Q[G]_L \isom \Q[x]/\Phi_r(x)\Q[x] \isom \Q(\bmu_r)$.  
The composition maps the maximal order $R_L$ isomorphically to 
the maximal order $\Z[\bmu_r]$, giving (ii).

Using (i) and (ii), there is a commutative diagram
$$
\xymatrix@C=15pt{
0 \ar[r] & \Phi_r(x)\Z[x]/(x^r-1)\Z[x] \ar[r] & \Z[x]/(x^r-1)\Z[x] 
    \ar[r] & \Z[x]/\Phi_r(x)\Z[x] \ar[r] & 0 \\
0 \ar[r] & \I_L^\perp \ar[r] \ar^\cong_{\eta}[u] & \Z[G] 
    \ar^\lambda[r] \ar^\cong_{\eta}[u] & R_L \ar[r] \ar^\cong[u] & 0
}
$$
with vertical isomorphisms, where the 
map $\lambda$ is induced by $\Q[G] \onto \Q[G]_L$.  
Since the top row is exact, so is the bottom row, giving (iii).

Let $\mu$ denote the M\"obius function. Then
\begin{equation}
\label{mi}
\Psi_r(x) = (x^r-1)/\Phi_r(x) = \prod_{d \mid r, d \ne 1} (x^{r/d}-1)^{-\mu(d)}.
\end{equation}
In $\Z[x]/\Phi_r(x)\Z[x]$, $x$ is a primitive $r$-th root of unity, 
so $x^{r/d}$ has order $d$, so $x^{r/d}-1$ is a unit in $\Z[x]/\Phi_r(x)\Z[x]$ 
unless $d$ is a prime power.  
When $d \ne 1$ is a prime power, $\mu(d) = -1$ if $d$ is prime, and 
$\mu(d) = 0$ otherwise.  
By (i), $\eta(\I_L)$ is generated by 
$\Psi_r(x)$, so by \eqref{mi}, the ideal $\lambda(\I_L)$ of $R_L$ is generated by 
$\prod_{\ell\mid r}(\zeta_\ell - 1)$.  
Since $\I_L \subseteq \Q[G]_L$, $\lambda$ is the identity map on $\I_L$.
This proves (iv).
\end{proof}  

\begin{thm}
\label{abprim0}
Suppose $L/\k$ is a cyclic extension of degree $r$, and $V$ is a 
commutative algebraic group over $\k$.  Then:
\begin{enumerate}
\item
$V_L$ is a commutative algebraic group of dimension $\varphi(r) \dim(V)$.
\item
If $V$ is connected then $V_L$ is connected.
\item
$V_L$ is isomorphic over $L$ to $V^{\varphi(r)}$.
\item
There is an injective ring homomorphism 
$R_L \otimes \End_\k(V) \hookto \End_\k(V_L)$.
\end{enumerate}
\end{thm}

\begin{proof}
Parts (i), (ii), and (iii) follow from Theorem \ref{easycor}, since 
$\rk_\Z \,\I_L = \dim \rho_L = \varphi(r)$.  Part (iv) follows from 
Proposition \ref{altthm}(iii) and Lemma \ref{newlem}(iii).
\end{proof}

\begin{lem}
\label{cycpollem}
The ideal of $\Z[x]/(x^r-1)\Z[x]$
generated by $\Phi_r(x)$ is also generated by each of the following sets
\begin{enumerate}
\item
$\{(x^r-1)/(x^d-1) : d \mid r, d \ne r\}$, 
\item
$\{(x^r-1)/(x^{r/\ell}-1) : \text{$\ell \mid r$, $\ell$ prime}\}$. 
\end{enumerate}
\end{lem}

\begin{proof}
The identity $x^r-1 = \prod_{d \mid r} \Phi_d(x)$ shows that 
$\Phi_r(x)$ divides $(x^r-1)/(x^d-1)$ for every divisor $d < r$ of $r$.  
On the other hand, Theorem 1 of \cite{deB} or  \cite{redei} shows that 
$\Phi_r(x)$ is a $\Z[x]$-linear combination of $\{(x^r-1)/(x^d-1) : d \mid r, d \ne r\}$. 
This proves (i).  Every element in the set (i) is divisible 
by one of the elements in its subset (ii), so this completes the proof.
\end{proof}

\begin{lem}
\label{cycpollem2}
Each of the sets $\{\Phi_r(\tau)\}$, $\{N_{L/F} : F \in \Omega_L\}$, 
$\{N_{L/F} : F \in \Omega'_L\}$ generates 
the ideal $\I_L^\perp \subseteq \Z[G]$.
\end{lem}

\begin{proof}
If $\k \subseteq F \subseteq L$ and $[F:\k] = d$, 
then 
$N_{L/F}$ goes to $(x^r-1)/(x^d-1)$
under the isomorphism $\Z[G] \isom \Z[x]/(x^r-1)\Z[x]$.  
Thus by Lemma \ref{cycpollem}, the three sets of this lemma generate
the same ideal of $\Z[G]$.
By Lemma \ref{newlem}(i), this ideal is $\I_L^\perp$.
\end{proof}

Recall that if $\alpha \in \Z[G]$, then $\alpha_V \in \End_\k(\Res^L_\k V)$ 
denotes its image under \eqref{rescpr2}.

\begin{thm}
\label{abprim}
Suppose $L/\k$ is a cyclic extension of degree $r$, and $V$ is a 
commutative algebraic group over $\k$.  Then:
\begin{enumerate}
\item
$V_L = \cap_{F \in \Omega_L} \ker(N_{L/F,V}) = \cap_{F \in \Omega'_L} \ker(N_{L/F,V}) 
    = \ker(\Phi_r(\tau)_V) \subseteq \Res^L_\k V$, 
where $\tau$ is any generator of $\Gal(L/\k)$.
\item
If $A$ is a commutative $\k$-algebra, then 
$$
V_L(A) \cong \{\alpha \in V(A \otimes_\k L) : 
    \text{$N_{L/F}(\alpha) = 0$ for every $F \in \Omega_{L}$}\}.
$$
In particular,
$$
V_L(\k) \cong \{\alpha \in V(L) : 
    \text{$N_{L/F}(\alpha) = 0$ for every $F \in \Omega_{L}$}\}.
$$
Both assertions also hold with $\Omega_L$ replaced by $\Omega'_L$.
\end{enumerate}
\end{thm}

\begin{proof}
Assertion (i) (resp., (ii)) follows from 
Lemma \ref{cycpollem2} and Proposition \ref{altthm}(i) (resp., (ii)).
\end{proof}

\begin{thm}
\label{abprim3}
Suppose $L/\k$ is a cyclic extension of degree $r$, and $V$ is a 
commutative algebraic group over $\k$.
Suppose that $\ell$ is prime and $g \in G_\k$.  Let $d$ 
be the order of the restriction of $g$ to $G:=\Gal(L/\k)$.  
If the characteristic polynomial of $g$ acting on $T_\ell(V)$ is 
$\prod_i (X-\alpha_i)$ with $\alpha_i \in \bar{\Q}_\ell$, 
then the characteristic polynomial of $g$ acting on $T_\ell(V_L)$ is 
$$
\prod_{i,\zeta}(X-\alpha_i \zeta)^{\varphi(r)/\varphi(d)}
$$
where $\zeta$ runs through all primitive $d$-th roots of unity.
\end{thm}

\begin{proof}
By Lemma \ref{newlem}(ii), 
the eigenvalues of the generator $\tau\in G$ acting on 
$\I_L \otimes\Q = R_L \otimes \Q$ are exactly the 
primitive $r$-th roots of unity in $\bar{\Q}$, each with multiplicity one.  
It follows that the eigenvalues of 
$g$ acting on $\I_L$ are the primitive $d$-th roots of unity, each 
with multiplicity $\varphi(r)/\varphi(d)$.  The result now follows from 
the isomorphism $T_\ell(V_L) \cong \I_L \otimes T_\ell(V)$
of Theorem \ref{vvntl}(iii).
\end{proof}

\begin{prop}
Suppose $L/\k$ is cyclic, $F$ and $M$ are extensions of $\k$ in $L$, 
$F \cap M = \k$, and $L = FM$.  
If $V$ is a commutative algebraic group over $\k$,
then $(V_{F})_{M} \cong V_L$ over $\k$.
\end{prop}

\begin{proof}
Let $d = [F:\k]$ and $e = [M:\k]$.  Then $de = r$. Since $L/\k$ 
is cyclic, $d$ and $e$ are relatively prime.  
By Lemma \ref{newlem}(ii,iv), there are isomorphisms of $\Z[G]$-modules 
$\I_F \cong \Z[\bmu_d]$, $\I_M \cong \Z[\bmu_e]$,  
and $\I_L \cong \Z[\bmu_r]$,
where the chosen generator $\tau$ of $G$ acts on the right-hand sides as 
multiplication by $\zeta_d$, $\zeta_e$, and $\zeta_r$, respectively,
and where the roots of unity are chosen so that $\zeta_d\zeta_e = \zeta_r$.
Then the natural map
$\Z[\bmu_d] \otimes_\Z \Z[\bmu_e] \isom \Z[\bmu_r]$
is an isomorphism of $\Z[G]$-modules.  Hence $\I_L \cong \I_M \otimes_\Z \I_F$,  
and the proposition follows from Proposition \ref{repeat}.
\end{proof}

\begin{rem}
Suppose $\k \subseteq F \subseteq L$.  
Let $N_{L/F} : \Z[G] \to \Z[G]$ denote multiplication by
$\sum_{h\in\Gal(L/F)} h$.
Then $N_{L/F}$ factors as 
$$
\xymatrix@C=40pt{
N_{L/F} : \Z[G] \ar@{->>}^-{R_{L/F}}[r] & \Z[\Gal(F/\k)]\, 
    \ar@{^(->}^-{\iota_{L/F}}[r] & \Z[G]
}
$$
where $R_{L/F}$ is the natural projection map.
Since $\ker(R_{L/F})$ and $\coker(\iota_{L/F})$ are torsion-free, 
Lemma \ref{twex} shows that the induced maps $R_{L/F,V}$ and $\iota_{L/F,V}$
in the composition
\begin{equation}
\label{triangle}
\xymatrix@C=40pt{
N_{L/F,V} : \Res^L_\k V \ar@{->>}^-{R_{L/F,V}}[r] & \Res^F_\k V\, 
    \ar@{^(->}^-{\iota_{L/F,V}}[r] & \Res^L_\k V
}
\end{equation}
are surjective and injective, respectively. 
In \cite{ssavs,hwpaper,ants}, 
the primitive subgroup of $\Res^L_\k V$ corresponding 
to $L$ was defined to be 
$T_L:=\cap_{\k \subseteq F \subsetneq L} \ker(R_{L/F,V})$.  
By \eqref{triangle}, $\ker(R_{L/F,V}) = \ker(N_{L/F,V})$. So 
when $L/\k$ is cyclic, $T_L=V_L=W_L$  (the last equality by 
Theorem \ref{abprim}(i)).
\end{rem}

\section{Semidirect products}
\label{sdprodsect}

Suppose for this section that $L/\k$ is a finite Galois extension, 
and $G := \Gal(L/\k)$ is a semidirect product $\Gamma \semidirect H$ 
of a normal cyclic subgroup $\Gamma = \Gal(L/K)$ of order $r$ by a 
subgroup $H = \Gal(L/M)$.  
There is a diagram
$$
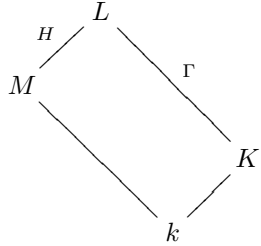
\xymatrix@R=15pt@C=15pt{
& L \\
M \ar@{-}^H[ur] \\
&&& K \ar@{-}_\Gamma[uull] \\
&& \k \ar@{-}[ur] \ar@{-}[uull]
}
$$
and we view $\Z[\Gamma]$ and $\Z[H]$ as subrings of $\Z[G]$, 
so $\Z[G] = \Z[\Gamma]\Z[H] = \Z[H]\Z[\Gamma]$.
Let $\rho_{L/K}$ be the (unique) irreducible faithful rational 
representation of $\Gamma$, and $\Q[\Gamma]_{L/K}$ the 
$\rho_{L/K}$-isotypic component of $\Q[\Gamma]$.  
Let $\I_L = \I_{L/K} \subseteq \Z[\Gamma]$ be the ideal 
$\Q[\Gamma]_{L/K} \cap \Z[\Gamma]$ of Definition \ref{vldef}, 
so $\I_L \in \mod{\Z}{K}$.  

In this section we will show 
(Theorem \ref{sdthm} below) that the commutative algebraic group 
$V_{L/K} = \I_L \otimes V \in \cag{\Z}{K}$ of Definition \ref{vldef} 
has a model over $\k$ of the form $\J_L \otimes V \in \cag{\Z}{\k}$ for a suitable 
right ideal $\J_L$ of $\Z[G]$.  This is needed 
for the applications in \cite{bigselmer}, in the case where 
$G$ is a dihedral group of order $2r$.

Define
$$
N_H := \sum_{h \in H} h \in \Z[H] \subseteq \Z[G].
$$

\begin{lem}
\label{sri}
The abelian group $\J_L := N_H \I_L$ is a saturated right ideal of $\Z[G]$.
\end{lem}

\begin{proof}
For $h \in H$, 
the representation $\rho_{L/K}^h$ of $\Gamma$ defined by 
$\rho_{L/K}^h(\gamma) = \rho_{L/K}(h \gamma h^{-1})$ is an irreducible 
faithful rational representation of $\Gamma$, so 
$\rho_{L/K}^h \cong \rho_{L/K}$.  
Hence $h \I_L h^{-1} = \I_L$ in $\Q[G]$, so for $h \in H$ 
and $\gamma \in \Gamma$ we have 
$$
N_H \I_L h = N_H h \I_L = N_H \I_L, \quad N_H \I_L \gamma = N_H \I_L,
$$
so $N_H \I_L \Z[G] = N_H \I_L$.  Since $\I_L \subseteq \Z[\Gamma]$ 
is saturated, $\J_L \subseteq \Z[G]$ is saturated.
\end{proof}

\begin{defn}
\label{ditwdef}
Define $V_{L/\k} := \J_L \otimes V$ where $\J_L := N_H \I_L$ as in Lemma \ref{sri}.
This definition depends on the subgroup $H$ of $G$; if necessary 
we will denote $V_{L/\k}$ by $V_{L/\k,H}$.  Theorem \ref{sdthm}
below shows that if $H'$ is another subgroup with $G = \Gamma \semidirect H'$, 
then $V_{L/\k,H'}$ is isomorphic to $V_{L/\k,H}$ over $K$.
\end{defn}

\begin{thm}
\label{sdthm}
Over $K$ there is an isomorphism $V_{L/\k} \cong V_{L/K}$, where 
$V_{L/K}$ (resp., $V_{L/\k}$) is given by Definition  \ref{vldef}
(resp., Definition \ref{ditwdef}).
\end{thm}

\begin{proof}
Left multiplication by $N_H$ is an isomorphism $\I_L \to \J_L$ of right $G_K$-modules.  
By Corollary \ref{obvcor} with $F = K$, 
$V_{L/K} = V_{\I_L}$ is isomorphic over $K$ to $V_{L/\k} = V_{\J_L}$.
\end{proof}

\section{Finite group actions on $\I \otimes V$}
\label{gpactsect}

In this section we study the action of symmetric groups on
the group varieties $\I \otimes V$.
When $V$ is an algebraic torus, these
results provide insights into some known cryptosystems 
(see \cite{hwpaper,ants}).

Fix a finite Galois extension 
$L/\k$ and let $G := \Gal(L/\k)$ (and $\O = \Z$).
Fix also a commutative algebraic group $V$ over $\k$ that is not isogenous 
to the trivial group, i.e., so that the natural map $\Z \to \End_\k(V)$ is injective.  
If $\sigma\in \End_\Z(\Z[G])$, 
let $\sigma_V \in \End_L(\Res^L_\k V)$ 
denote the endomorphism given by the functorial map 
of Corollary \ref{jdrem}(i) (with $\I = \J = \Z[G]$).
If $\I$ is a saturated right ideal of $\Z[G]$, view 
$\I \otimes V \subseteq \Res^L_\k V$ via Lemma \ref{twex} (with
$\J= \Z[G]$) and Proposition \ref{rescpr1}.  

\begin{lem}
\label{preboxes}
Suppose that $\I$ is a saturated right ideal of $\Z[G]$, and 
$\sigma \in \End_\Z(\Z[G])$.  
Then the following are equivalent:
\begin{enumerate}
\item
$\sigma(\I) \subseteq \I$.
\item
${\sigma}_V(\I \otimes V) \subseteq \I \otimes V$.
\end{enumerate}
\end{lem}

\begin{proof}
If $\sigma(\I) \subseteq \I$ then $\sigma|_\I\in\End_\Z(\I)$. 
By the functoriality of $\I \mapsto \I \otimes V$, 
we have $\sigma_V|_{\I \otimes V}\in\End_L(\I \otimes V)$.  
Thus (i) $\Rightarrow$ (ii).

Conversely, suppose ${\sigma}_V(\I \otimes V) \subseteq \I \otimes V$ and
let $\lambda : \Z[G] \to \End_\Z(\Z[G])$ denote the  
map that sends $\alpha \in \Z[G]$ to 
left multiplication by $\alpha$. 
Suppose $\alpha \in \I$ and $\beta \in \I^\perp$.  
Then $\alpha \Z[G] \subseteq \I$, so 
$\alpha_V \in \End_\k(\Res^L_\k V)$ factors through $\I \otimes V$, i.e., 
$\alpha_V(\Res^L_\k V) \subseteq \I \otimes V$. 
Therefore
$$
(\lambda(\beta) \circ \sigma \circ \lambda(\alpha))_V(\Res^L_\k V) 
    = \beta_V(\sigma_V(\alpha_V(\Res^L_\k V))) 
    \subseteq \beta_V(\I \otimes V) = 0
$$
by Proposition \ref{altthm}(i).  By Corollary \ref{jdrem}(iii), 
the map $\End_\Z(\Z[G]) \to \End_{\ks}(\Res^L_\k V)$ 
is injective, so 
$\lambda(\beta) \circ \sigma \circ \lambda(\alpha) = 0$, 
and thus
$
\beta \cdot \sigma(\alpha) 
    = (\lambda(\beta) \circ \sigma \circ \lambda(\alpha))(1) = 0.
$
Therefore $\I^\perp \sigma(\I) = 0$, so $\sigma(\I) \subseteq \I$ by 
Proposition \ref{endl}(ii).  Thus (ii) $\Rightarrow$ (i).
\end{proof}

Let $\Sigma_H$ denote the group of permutations of a set $H$. 
If $\sigma \in \Sigma_G$, let 
$\hat{\sigma} \in \Aut_\Z(\Z[G])$ denote the automorphism induced by $\sigma$,
and let $\hat{\sigma}_V  \in \Aut_L(\Res^L_\k V)$ denote  the corresponding 
automorphism of $\Res^L_\k V$ .

\begin{lem}
\label{sigautprop}
Suppose that $L/\k$ is cyclic and $\sigma \in \Sigma_G$.  
Then the restriction of $\hat{\sigma}_V$ to $V_L$ 
is an automorphism of $V_L$ if and only if 
\begin{itemize}
\item[(*)]
for every $g \in G$ and subgroup $H \subseteq G$ of prime order, 
$\sigma(gH) = \sigma(g)H$.
\end{itemize}
\end{lem}

\begin{proof}
Since $\sigma$ has finite order, the restriction of $\hat{\sigma}_V$ 
to $V_L$ is an automorphism 
if and only if $\hat{\sigma}_V(V_L) \subseteq V_L$,  
which by Lemma \ref{preboxes} happens  
if and only if $\hat{\sigma}(\I_L) \subseteq \I_L$.
Write $G = G_1 \times \cdots \times G_t$ where each
$G_i$ is of prime power order, $|G_i|=p_i^{r_i}$,  ordered so that
$p_1 < \ldots < p_t$.  For  $1\le i\le t$, 
let $H_i$ be the subgroup of $G_i$ of order $p_i$, and let  
$N_{H_i} = \sum_{h \in H_i} h \in \Z[G]$.
By Lemma \ref{cycpollem2}, 
\begin{equation}
\label{iomegaeq}
\I_L^\perp = \sum_{i=1}^t \Z[G]N_{H_i}.
\end{equation}
If $\sigma$ satisfies (*) then 
$\hat{\sigma}(N_{H_i} \alpha) = N_{H_i} \cdot \hat{\sigma}(\alpha)$
for every $\alpha \in \Z[G]$ and every $i$, so  
$$
\I_L^\perp\cdot\hat{\sigma}(\I_L) =
\sum_{i}\Z[G]N_{H_i} \cdot \hat{\sigma}(\I_L) 
    = \sum_{i}\Z[G]\cdot\hat{\sigma}(N_{H_i} \cdot \I_L) = 0,
$$
since $N_{H_i} \cdot \I_L = 0$ for all $i$.
By Proposition \ref{endl}(ii) we conclude that 
$\hat{\sigma}(\I_L) \subseteq \I_L$, so by Lemma \ref{preboxes}, 
$\hat{\sigma}_V |_{V_L}\in\Aut(V_L)$.

Conversely, suppose (*) fails to hold for some $H$.  
Take $j$ minimal so that there is a $\gamma \in G$ with 
$\sigma(\gamma H_j) \ne \sigma(\gamma)H_j$.  
Replacing $\sigma$ by $\tau_{\sigma(\gamma)^{-1}} \circ \sigma \circ \tau_\gamma$ 
(where $\tau_g \in \Sigma_G$ is left multiplication by $g \in G$)
we may 
assume without loss of generality that $\sigma(1) = 1$, $\sigma(H_j) \ne H_j$, 
and $\sigma(gH_i) = \sigma(g)H_i$ for all $g\in G$ and $i<j$.
It follows that
$\sigma^{-1}(gH_i) = \sigma^{-1}(g)H_i$ for all $g\in G$ and $i<j$, so
\begin{equation}
\label{dc}
\sigma^{-1}(g\prod_{i<j}H_i) = \sigma^{-1}(g)\prod_{i<j}H_i \quad \text{for every $g \in G$}.
\end{equation}

Let $\pi_i: G \to G_i$ be the projection map.
For $1\le i\le t$, fix $1 \neq \delta_i\in H_i$
such that
$$
\delta_i \notin
\begin{cases}
\pi_i(\sigma^{-1}(H_j)) & \text{ if $i>j$}, \\
\sigma^{-1}(H_j) & \text{ if $i=j$}, \\
\pi_i(\sigma^{-1}(H_j)\cap\delta_j\prod_{i<j}H_i) & \text{ if $1<i<j$}.
\end{cases}
$$
The first is possible since $p_i > p_j$ if $i > j$; the second 
since $\sigma(H_j) \ne H_j$; and the third because it follows from
\eqref{dc} that the elements of $\sigma^{-1}(H_j)$ lie in distinct cosets 
of $\prod_{i<j}H_i$, and $|H_i|=p_i \ge 3$ if $i > 1$.

If $S \subseteq \{1,\cdots,t\}$, let $\delta_S=\prod_{i\in S}\delta_i$.
Note that $\delta_S=1$ if and only if $S=\emptyset$. 
We claim that if $\sigma(\delta_S)\in H_j$, then either
$S=\emptyset$, or else
$S=\{1,j\}$ and $j\neq 1$ (in which case $\delta_S=\delta_1\delta_j$).
To prove the claim, suppose $S\neq\emptyset$ (so $\delta_S \ne 1$) 
and $\sigma(\delta_S)\in H_j$.
Then $\pi_i(\delta_S)\in \pi_i(\sigma^{-1}(H_j))$.
Note that $\pi_i(\delta_S)$ is $\delta_i$ if $i\in S$ and is $1$ otherwise.  
By our constraints on the $\delta_i$, if $i > j$ then 
$i \notin S$. 
If $j \notin S$, then applying \eqref{dc} with $g=1$ gives
$\sigma(\delta_S)\in (\prod_{i<j}H_i) \cap H_j =\{1\}$,
contradicting that $\sigma(1)=1$ and $\delta_S\neq 1$.
Thus $j \in S$, 
and again by our constraints,  if $1 < i < j$ then $i \notin S$.  
Since $\sigma(\delta_j)\notin H_j$, we cannot have $S=\{j\}$.
We have thus proved the claim.

Let 
$\alpha := \prod_{i=1}^t(1-\delta_i) 
    = \sum_S (-1)^{|S|}\delta_S \in \Z[G]$,
with $S$ running over subsets of $\{1,\cdots,t\}$.
By the claim above,  
$\hat{\sigma}(\alpha) = \sum_S (-1)^{|S|}\sigma(\delta_S) = 
1 + \sum_{g\notin H_j}a_g g$ or
$1 + \sigma(\delta_1\delta_j) + \sum_{g\notin H_j}a_g g$
with $a_g \in \Z$.
It follows that $N_{H_j} \cdot \hat{\sigma}(\alpha) \ne 0$,
since the coefficient of the identity element is either $1$ or $2$, 
so $\hat{\sigma}(\alpha) \notin \I_L$ by Lemma \ref{cycpollem2}.
Since $\delta_i \in H_i$, we have $N_{H_i}(1-\delta_i)=0$ for all $i$.
Thus by \eqref{iomegaeq}, $\I_L^\perp \alpha = 0$, so 
by Proposition \ref{endl}(ii), $\alpha \in \I_L$.  
Therefore $\hat\sigma(\I_L) \not\subseteq \I_L$, so by 
Lemma \ref{preboxes}, $\hat{\sigma}_V(V_L) \not\subseteq V_L$.
\end{proof}

If $|G|$ is squarefree, and $H$ is a subgroup of $G$, 
then there is a unique subgroup 
$J \subseteq G$ such that $G=H\times J$, and
this decomposition induces an inclusion $\Sig{H}\subseteq\Sig{G}$.

\begin{thm}
Suppose $L/\k$ is cyclic of squarefree degree,
$|G|=p_1\cdots p_t$ with distinct primes $p_i$, 
$H_i$ is the subgroup of $G$ of order $p_i$,
and $\sigma \in \Sigma_G$.  
Then $\hat{\sigma}_V|_{V_L}\in\Aut(V_L)$ if and only if 
$\sigma \in \prod_{i=1}^t \Sig{H_i} 
(\subseteq 
\Sig{G})$.
\end{thm}

\begin{proof}
Suppose $\sigma \in \prod_i \Sig{H_i}$. It is easy to see that
for every $g\in G$ and every $i$, $\sigma(gH_i)=\sigma(g)H_i$. 
By Lemma \ref{sigautprop}, $\hat{\sigma}_V|_{V_L}\in\Aut(V_L)$.

Conversely, suppose $\hat{\sigma}_V|_{V_L}\in\Aut(V_L)$.
By Lemma \ref{sigautprop},
 $\sigma(gH_i)=\sigma(g)H_i$ 
for all $g\in G$ and all $i$.
Let $\pi_i : G \to H_i$ denote the projection, and let
$\sigma_i=\sigma|_{H_i} : H_i \to G$.
Let $\tau_i = \pi_i\circ\sigma_i\in\Sig{H_i}$.
It follows easily that 
$\sigma=\prod_{i=1}^t \tau_i\circ\pi_i\in\sigma\in\prod_i \Sig{H_i}$.
\end{proof}

\appendix
\section*{Appendix.  More general construction}

Although in the above discussion we restrict to the case of free $\O$-modules 
$\I$, the tensor product construction  $(\I, V) \mapsto  \ItV$ 
in the appropriate category of sheaves for the {\'e}tale topology 
(as alluded to in the introduction) is quite general. Moreover, 
this more general construction can also be formulated in fairly 
concrete terms. 
For example, suppose that $\O$ is a commutative noetherian ring, 
$V \in \cagok$, and $\I$ is a finitely generated $\O$-module with 
a continuous right action of $G_\k$, but do not assume that 
$\I$ is a free $\O$-module.
Let $L$ be a finite Galois extension of $\k$ such that $G_L$ 
acts trivially on $\I$, and let $G := \Gal(L/\k)$.  
Since $\O$ is noetherian, there is   
an $\O[G]$-presentation of $\I$, i.e., an exact sequence
$$
\O[G]^a \map{~\psi~} \O[G]^b \too \I \too 0
$$ 
of right $\O[G]$-modules.  By basic properties of the functor $V \mapsto \Res^L_\k V$ 
(or for example, by Corollary \ref{jdrem}(ii) 
and Proposition \ref{rescpr1}), $\psi$ induces a $\k$-homomorphism 
$$
\psi_V : (\Res^L_\k V)^a \to (\Res^L_\k V)^b,
$$ 
and we can define
$$
\I \otimes_\O V := \coker(\psi_V) \in \cagok.
$$
One can show that this definition is independent of the choice 
of $L$ and of the presentation of $\I$, and it agrees with Definition \ref{ItVdef} 
if $\I$ is a free $\O$-module.  Further, 
without the assumption that the $\O$-modules are free, 
Theorem \ref{ithm}, Corollaries \ref{obvcor}, \ref{cmcor}, and \ref{decompcor}, 
and Lemma \ref{isoglem} all remain true verbatim, 
Proposition \ref{newjdiag} and Corollary \ref{jdrem} hold 
if $\I$ and $\J$ are projective $\O$-modules, 
Theorem \ref{vvntl} holds if $\I$ is a projective $\O$-module, and Lemma \ref{twex} 
holds if $\I/\J$ is a projective $\O$-module.

   This definition of $\ItV$ is essentially the same
   as Conrad's definition of his $\I \otimes_{\O[G]} \Res^L_\k V$
   in Theorem 7.2 of \cite{conrad}, using the action of $\O[G]$ on
   $\Res^L_\k V$ given by \eqref{rescpr2} above.

\end{document}